%% file: main.tex
\documentclass[final,reqno]{amsart}

\input{macro}

\begin{document}
\title{Connection-twisted modules and logarithmic twisted modules of vertex algebras}
\author{Takumi Iwane}
\date{\today}
\address{Graduate School of Mathematics, Nagoya University.
 Furocho, Chikusaku, Nagoya, Japan, 464-8602.}
\email{takumi.iwane.c8@math.nagoya-u.ac.jp}
\begin{abstract}
    We formulate the notion of connection-twisted modules for $\mathbb{C}$-graded vertex algebras, as suggested by Feigin and Lentner, and investigate their fundamental properties. Let $G$ be an algebraic group acting on a vertex algebra $V$ by automorphisms. We establish an equivalence between the category of $\nabla$-twisted $V$-modules and the category of $g$-twisted $V$-modules, where $\nabla$ is a regular singular $\mathfrak{g}$-connection and $g\in G$ is its monodromy. For symplectic fermion and affine vertex algebras, we explicitly compute the commutation relations and the action of the Virasoro algebra on connection-twisted modules, recovering the formulas obtained by Feigin and Lentner \cite{FL24,FL25}.    
\end{abstract}

\maketitle


\section{Introduction}

Vertex algebras provide an algebraic framework for describing the chiral sector of two-dimensional conformal field theories and have revealed deep connections among representation theory, quantum groups, algebraic geometry, and mathematical physics (standard references \cite{Kac98,LL04,FBZ04}). Under suitable assumptions, the category of modules over a vertex algebra admits the structure of a braided tensor category (see \cite{HL13}) and is closely related to the representation theory of quantum groups. A prominent example is the Kazhdan--Lusztig correspondence \cite{KL93,KL94}, which relates representations of affine vertex algebras to those of quantum groups. In logarithmic conformal field theory, the logarithmic Kazhdan--Lusztig correspondence predicts a close relationship between representation categories of certain logarithmic vertex algebras and quantum groups at roots of unity \cite{FGST06,CLR21,CLR23}. In particular, the representation category of the Feigin–Tipunin algebra is conjectured to be equivalent to that of the corresponding small quantum group. This conjecture is one of the central problems in the study of non-semisimple braided tensor categories.

In recent years, vertex algebras arising from the quantum geometric Langlands correspondence have attracted considerable attention as a framework for understanding such relations in a broader context. Motivated by the gauge-theoretic interpretation of geometric Langlands duality due to Kapustin and Witten \cite{KW07} and subsequent developments \cite{GR19,FG20}, certain vertex algebras arise as chiral algebras associated with two-dimensional junctions between boundary conditions in topologically twisted four-dimensional gauge theories. These vertex algebras, known as quantum geometric Langlands kernels, play an important role in the study of the quantum geometric Langlands correspondence and dualities of vertex algebras.

In \cite{FL24}, the authors investigated the semiclassical limits of generalized quantum Langlands kernels and their representation theory parametrized by spaces of connections. In general, a family of vertex algebras may acquire a large center at a special value of its deformation parameter. A classical example is the Feigin--Frenkel center of an affine vertex algebra at the critical level \cite{FF92,Fre05,Fre07}. This center is identified with the algebra of functions on the space of opers associated with the Langlands dual Lie algebra, thereby connecting the representation theory of vertex algebras with the geometry of differential equations. Similarly, Feigin and Lentner showed that the large-level limits of generalized quantum Langlands kernels acquire large centers containing the algebra of functions on the space $\mathrm{Conn}_{\mathfrak g}(\mathbb{D})$ of $\mathfrak g$-valued regular connections on the disk $\mathbb{D}$.

More generally, let $\mathcal V$ be a vertex algebra whose center contains the commutative algebra $\shO(\mathrm{Conn}_{\mathfrak g}(\mathbb{D}))$ of functions on $\mathrm{Conn}_\frg(\mathbb{D})$. It is then natural to study the representation theory of $\mathcal V$ fiberwise over $\mathrm{Conn}_{\mathfrak g}(\mathbb{D})$, with each fiber corresponding to a $\mathfrak g$-valued connection $\nabla$. In the setting considered by \cite{FL24,FL26}, the fiber over the zero connection is expected to recover a vertex algebra $V$, while the other fibers are expected to be related to the representation theory of $V$. The authors proposed that these fibers should be understood not as ordinary $V$-modules, but rather as modules twisted by the corresponding connections.

This perspective leads naturally to the theory of logarithmic twisted modules over vertex algebras (developed in \cite{Hua10,Bak16}). Let $V$ be a vertex algebra equipped with an action of an algebraic group $G$ by automorphisms. For each $g\in G$, one can consider logarithmic $g$-twisted $V$-modules, which describe twisted sectors in which fields acquire a nontrivial $g$-action upon analytic continuation around the origin. Such twisted sectors play a fundamental role in orbifold conformal field theory and, under suitable assumptions, form a $G$-crossed braided tensor category (see \cite{McR21,DH26,Xu26}). Instead of specifying an automorphism $g$, one may consider a connection
\[
    \nabla=d+A(z)dz
\]
and seek a notion of modules twisted by the connection itself. For regular singular connections, Feigin and Lentner suggested that such modules should be equivalent to logarithmic twisted modules through the monodromy of the connection. In contrast, irregular singular connections are expected to give rise to a richer representation theory that cannot be fully captured by monodromy alone.


These developments suggest the need for a general algebraic framework for vertex algebra modules twisted by connections. In particular, two fundamental questions arise: how gauge transformations relate the module categories associated with different connections, and how modules twisted by regular singular connections are related to the existing theory of logarithmic twisted modules. The purpose of this paper is to develop such a framework and answer these questions.

We introduce the notion of a $\nabla$-twisted module $M=(M,Y_M^\nabla)$ for a $\mathfrak g$-valued connection $\nabla=d+A(z)dz$ (see \cref{dfn:conn tw mod}). The key ingredient of our formulation is the formal parallel transport operator:
\[
    P_\nabla(w,z)=e^{z(\partial_w+A(w))}1=1+zA(w)+z^2\frac{A'(w)+A(w)^2}{2}+\cdots
\]
which allows us to incorporate the connection directly into the Jacobi identity. More precisely, we require the $\nabla$-twisted Jacobi identity:
\[
\begin{split}
&z_0^{-1}\delta\left(\frac{z_1-z_2}{z_0}\right)
Y_M^\nabla(a,z_1)Y_M^\nabla(b,z_2)+(-1)^{|a||b|}(-z_0)^{-1}
\delta\left(\frac{z_2-z_1}{-z_0}\right)
Y_M^\nabla(b,z_2)Y_M^\nabla(a,z_1)\\
&=z_2^{-1}\delta\left(\frac{z_1-z_0}{z_2}\right)
Y_M^\nabla\bigl(Y(P_\nabla(z_2,z_0)a,z_0)b,z_2\bigr).
\end{split}
\]
This identity provides a natural deformation of the usual Jacobi identity and leads to connection-dependent analogues of the Borcherds identity, the $n$-th product formula, and the translation property.

Our first main result establishes the compatibility of this definition with gauge transformations.
\begin{mainthm}[{\cref{thm:gauge transf,cor:main A}}]
    Let $V$ be a $\bbC$-graded vertex algebra equipped with a vertex algebraic $G$-action, and let $\nabla$ be a $\mathfrak g$-valued connection. Denote by $\cMod^\nabla(V)$ the category of $\nabla$-twisted $V$-modules. Then, for any $g(z)\in G(\mathbb{C}\dpr{z})$, the gauge transformation induces an equivalence of categories
    \[
        (-)^g\colon\cMod^\nabla(V)\to\cMod^{\nabla^g}(V),
    \]
    where $\nabla^g=g(z)^{-1}\nabla g(z)$, with quasi-inverse $(-)^{g^{-1}}$
\end{mainthm} 
In particular, the category of connection-twisted modules depends only on the gauge-equivalence class of the connection, up to equivalence.

For regular singular connections, we extend the above correspondence to a fundamental solution involving fractional powers and logarithms. This leads to our second main result, which gives an affirmative answer to the question raised by Feigin and Lentner \cite{FL25,FL26}.

\begin{mainthm}[{\cref{thm:conn tw mod to log tw mod,thm:log tw mod to conn tw mod,cor:main B}}]
    Let $V$ be a $\bbC$-graded vertex algebra equipped with a vertex algebraic $G$-action, and let $\nabla=d+A(z)dz$ be a regular singular $\mathfrak{g}$-connection. Suppose that $F(z)\in G(\widetilde{\shK}_z)$ is a fundamental solution satisfying
    \[
        F(z)^{-1}\nabla F(z)=d,
    \]
    and let $g\in G(\mathbb{C})$ be its monodromy. Then the functor
    \[
        (-)^F\colon\cMod^\nabla(V)\lto\cMod^g(V)
    \]
    is an equivalence of categories, with quasi-inverse $(-)^{F^{-1}}$. Here, $\cMod^g(V)$ denotes the category of logarithmic $g$-twisted $V$-modules.
\end{mainthm}
Thus, for regular singular connections, the category of connection-twisted modules is determined, up to equivalence, by the monodromy automorphism. In particular, our formulation recovers the existing theory of logarithmic twisted modules in the regular singular case.

We further illustrate our theory through the examples of the symplectic fermion vertex superalgebra and affine vertex algebras. In the case of symplectic fermion, we verify that our construction reproduces the connection-dependent structures investigated in \cite{FL25}. For affine vertex algebras, we obtain the following description:
\begin{mainthm}[{\cref{prp:deform mode alg aff va}}]
    Let $V^k(\mathfrak{g})$ be an affine vertex algebra and let $\nabla=d+A(z)dz$ be a $\mathfrak{g}$-valued connection. Define the Lie algebra $\widehat{\mathfrak g}_\nabla=\mathfrak g((t))\oplus\mathbb CK$ with Lie bracket
\[
    [x_m^\nabla,y_n^\nabla]
    =[x,y]_{m+n}^\nabla
    +(A_{m+n},[x,y])K
    +m\delta_{m+n,0}(x,y)K,
\]
for $x_m^\nabla=xt^m$ and $K$ is central. Then the Fourier modes of $Y_M^\nabla(x,z)$ on any $\nabla$-twisted $V^k(\mathfrak{g})$-module satisfy the defining relations of $\widehat{\mathfrak g}_\nabla$ at level $k$.
\end{mainthm}

Since the additional term $(A_{m+n},[x,y])$ is a Lie algebra $2$-coboundary, we also show that $\widehat{\mathfrak g}_\nabla$ is isomorphic to the ordinary affine Lie algebra $\widehat{\mathfrak g}$ (\cref{prp:isom aff Lie alg}). In \cite{FL24}, the authors expect that for any $\frg$-connection $\nabla$, the category $\cMod^\nabla(V^k(\frg)$ is equivalent to the category of $V^k(\frg)$-modules. If we prove that $\cMod^\nabla(V^k(\frg))$ is equivalent to the category of continuous $\wh{\frg}_\nabla$-modules, then we argue that the above expectation is true. To describe $\cMod^\nabla(V^k(\frg))$ in the mode algebra language is ongoing work.

Furthermore, we determine the Virasoro action on $\nabla$-twisted $V^k(\mathfrak g)$-modules through a deformed Sugawara construction (\cref{prp:conf str on tw aff}):
\begin{equation}\label{eq:deformed Sugawara vector}
    L_{\mathrm{sug}}^\nabla(z)
    =
    \frac{\sum_i\nop{x_i(z)x^i(z)}}{2(k+h^\vee)}
    -(A(z))(z)
    +\frac{k h^\vee}{2}(A(z),A(z)).
\end{equation}
Here the notation is as in \cref{prp:isom aff Lie alg}. This formula is similar to that of \cite{FL24}
\footnote{The discrepancy between our formula and that of \cite{FL24} appears to arise from the minus sign in front of the Sugawara vector in the proof of \cite[Theorem~5.3]{FL24}.}.

\subsection*{Organization} In \cref{s:preliminaries}, we recall the basic definitions and properties of vertex algebras and their modules. In \cref{s:def and basic properties}, we introduce connection-twisted modules using the formal parallel transport operator and establish their fundamental properties, including the twisted Borcherds identity, the deformed $n$-th product formula, and the translation property. In \cref{s:gauge transformation}, we study gauge transformations and prove that gauge-equivalent connections give rise to equivalent categories of twisted modules. In \cref{s:regular singular case}, we consider regular singular connections and establish an equivalence between connection-twisted modules and logarithmic twisted modules associated with their monodromy. In \cref{s:example}, we apply our theory to the symplectic fermion vertex superalgebra and affine vertex algebras, recovering the corresponding deformed mode algebras and Virasoro actions. In Appendix~A, we use Tannakian formalism to establish the existence of group-valued fundamental solutions for regular singular connections.

\section{Preliminaries}\label{s:preliminaries}
In this section, we review the notions of vertex algebras and their modules, based on \cite{LL04} (see also \cite{Kac98,FBZ04}). Throughout this paper, we refer to vertex superalgebras simply as vertex algebras.

\subsection{Vertex algebras}
For a vector superspace $V$, let $V\dbr{z}$ and $V\dpr{z}$ denote the spaces of formal power series and formal Laurent series with coefficients in $V$, respectively. Note that $V\otimes\bbC\dbr{z}\subsetneq V\dbr{z}$ and $V\otimes\bbC\dpr{z}\subsetneq V\dpr{z}$ for an infinite-dimensional vector superspace \(V\).

To formulate the definition of vertex algebras, we introduce the \emph{formal delta function}: 
\[
    \delta(z)\ceq \sum_{n\in \bbZ}z^n.
\]    
The following two properties of the formal delta function are fundamental and useful: 
\begin{prp}\label{prp:basic property of formal delta funtion}
    Let $V$ be a vector superspace. For any $f(z_0,z_1,z_2)\in V\dbr{z_0,z_1,z_2}[z_0^{-1},z_1^{-1},z_2^{-1}]$, we have
    \[
        z_0^{-1}\delta\left(\frac{z_1-z_2}{z_0}\right)f(z_0,z_1,z_2)=z_0^{-1}\delta\left(\frac{z_1-z_2}{z_0}\right)f(z_1-z_2,z_1,z_2),
    \]
    where $f(z_1-z_2,z_1,z_2)$ on the right hand side is the expansion in $\abs{z_1}>\abs{z_2}$, i.\ e.\ 
    \[
        (z_1-z_2)^n=\sum_{k\geq 0}(-1)^k\binom{n}{k}z_1^{n-k}z_2^k.
    \]
\end{prp}

\begin{prp}
    The formal delta function satisfies the following identities:
    \begin{align}
        &z_2^{-1}\delta\left(\frac{z_1-z_0}{z_2}\right)=z_1^{-1}\delta\left(\frac{z_2+z_0}{z_1}\right) \label{eq: skew symmetry}\\
        &z_0^{-1}\delta\left(\frac{z_1-z_2}{z_0}\right)+(-z_0)^{-1}\delta\left(\frac{z_2-z_1}{-z_0}\right)=z_2^{-1}\delta\left(\frac{z_1-z_0}{z_2}\right) \label{eq: Jacobi identity}
    \end{align}
\end{prp}

Here, we introduce the notion of fields:
\begin{dfn}\label{dfn: local field and n-th product}
    For a vector superspace $V$, an element of $\Hom(V,V\dpr{z})$ is called a \emph{field} on $V$. If fields $a(z)$ and $b(w)$ satisfy that $(z-w)^N[a(z),b(w)]=0$ for $N\gg 0$, then $a(z)$ and $b(w)$ are \emph{local}. Then 
    \[
        a(w)_{(n)}b(w)\ceq \oint dz \bigl((z-w)^na(z)b(w)-(-1)^{\abs{a}\abs{b}}(w-z)^nb(w)a(z)\bigr)
    \]
    is again a field, called the \emph{$(n)$-th product}. Here, $\oint dz$ is the operation of taking the coefficient of $z^{-1}$. In particular, 
    \[
        \begin{split}
            &\nop{a(w)b(w)}=a(w)_{(-1)}b(w)=a(w)_+b(w)+(-1)^{\abs{a}\abs{b}}b(w)a(w)_-,\\
            &a(w)_+=\sum_{n<0}a_nz^{-n-1},\quad a(w)_-=\sum_{n\geq 0}a_nz^{-n-1}
        \end{split}
    \]
    is called the \emph{normally ordered product}.
\end{dfn}

The following is a fundamental property of $(n)$-th product:
\begin{prp}\label{prp:fundamental property of n-th product}
    For local two fields $a(z)$ and $b(w)$, $a(z)b(w)$ is expressed as
    \[
        a(z)b(w)=\sum_{n<N}(z-w)^{-n-1}a(w)_{(n)}b(w)
    \]
    for some $N\geq 0$, where the right hand side is expanded in $\abs{z}>\abs{w}$. 
\end{prp}

\begin{dfn}
    A \emph{vertex algebra} $V=(V,\one,Y)$ consists of the following data:
    \begin{itemize}
        \item 
        $V$ is a vector superspace,
        \item 
        $Y\colon V\to \Hom(V,V\dpr{z});\ a\mto Y(a,z)=\sum_{n\in\bbZ}a_{(n)}z^{-n-1}$ is a linear map, preserving parities, 
    \end{itemize}
    satisfying that
    \begin{enumerate}[series=axiomva]
        \item 
        $Y(\one,z)=\id_V$, $Y(a,z)\one\in V\dbr{z}$, and $Y(a,z)\one|_{z=0}=a$,
        \item 
        Jacobi identity:
        \[
        \begin{split}
            &z_0^{-1}\delta\left(\frac{z_1-z_2}{z_0}\right)Y(a,z_1)Y(b,z_2)+(-1)^{\abs{a}\abs{b}}(-z_0)^{-1}\delta\left(\frac{z_2-z_1}{-z_0}\right)Y(b,z_2)Y(a,z_1)\\
            &=z_2^{-1}\delta\left(\frac{z_1-z_0}{z_2}\right)Y(Y(a,z_0)b,z_2).
        \end{split}
        \]
    \end{enumerate}
\end{dfn}

For later use in \cref{s:def and basic properties}, we introduce the notion of $\bbC$-graded vertex algebras:
\begin{dfn}
    A \emph{$\bbC$-graded vertex algebra} is a vertex algebra $V$ equipped with a $\bbC$-grading $V=\bigoplus_{\Delta\in\bbC}V_\Delta$ satisfying the following properties:
    \begin{enumerate}[resume=axiomva]
        \item 
        for any $\Delta\in\bbC$, $V_{\Delta-n}=0$ for sufficiently large $n\in\bbZ$,
        \item
        $\one\in V_0$,
        \item
        for a homogeneous element $a\in V$, the grading of $a_{(n)}$ is $\Delta_a-n-1$, where $\Delta_a$ denotes the grading of $a$
    \end{enumerate}
\end{dfn}

The following properties are fundamental (see \cite[\S~3.3, 3.4]{LL04}).

\begin{prp}\label{prp: axiomva}
    The axiom (2) is equivalent to
    \begin{enumerate}[label=(2-\arabic*)]
        \item 
        weak commutativity (locality): for any $a,b\in V$, there exists $N\geq 0$ such that
        \[
            (z_1-z_2)^NY(a,z_1)Y(b,z_2)=(-1)^{\abs{a}\abs{b}}(z_1-z_2)^NY(b,z_2)Y(a,z_1),
        \]
        \item
        weak associativity: for any $a,b\in V$, there exists $N\geq 0$ such that
        \[
            (z_0+z_2)^NY(a,z_0+z_2)Y(b,z_2)=(z_0+z_2)^NY(Y(a,z_0)b,z_2).
        \]
    \end{enumerate}
\end{prp}

\begin{prp}\label{prp: trans equiv va}
    We define a linear operator $T$ on $V$ by
    \[
        Ta\ceq\oint dz\ z^{-2}Y(a,z)\one=a_{(-2)}\one.
    \]
    Then, we have
    \[
        [T,Y(a,z)]=Y(Ta,z)=\frac{d}{dz}Y(a,z),
        \quad
        a(z)\one=e^{zT}a
    \]
\end{prp}

\subsection{Vertex algebra modules}
\begin{dfn}
    A \emph{$V$-module} $M=(M,Y_M)$ consists of the following data:
    \begin{itemize}
        \item
        $M$ is a vector superspace,
        \item
        $Y_M\colon V\to \Hom(M,M\dpr{z})$ is a linear map, preserving parities,
    \end{itemize}
    \begin{enumerate}
        \item
        $Y_M(\one,z)=\id_M$
        \item
        Jacobi identity:
        \[
        \begin{split}
            &z_0^{-1}\delta\left(\frac{z_1-z_2}{z_0}\right)Y_M(a,z_1)Y_M(b,z_2)+(-1)^{\abs{a}\abs{b}}(-z_0)^{-1}\delta\left(\frac{z_2-z_1}{-z_0}\right)Y_M(b,z_2)Y_M(a,z_1)\\
            &=z_2^{-1}\delta\left(\frac{z_1-z_0}{z_2}\right)Y_M(Y(a,z_0)b,z_2).
        \end{split}
        \]
    \end{enumerate}
\end{dfn}

\begin{prp}\label{prp: axiommod}
    The axiom (2) is equivalent to
    \begin{enumerate}[label=(2-\arabic*)]
        \item 
        weak commutativity (locality): for any $a,b\in V$, there exists $N\geq 0$ such that
        \[
            (z_1-z_2)^NY_M(a,z_1)Y_M(b,z_2)=(-1)^{\abs{a}\abs{b}}(z_1-z_2)^NY_M(b,z_2)Y_M(a,z_1),
        \]
        \item
        weak associativity: for any $a,b\in V$, there exists $N\geq 0$ such that
        \[
            (z_0+z_2)^NY_M(a,z_0+z_2)Y_M(b,z_2)=(z_0+z_2)^NY_M(Y(a,z_0)b,z_2).
        \]
    \end{enumerate}
\end{prp}

\begin{prp}\label{prp: trans equiv mod}
    For a $V$-module $M$, we have
    \[
        Y_M(Ta,z)=\frac{d}{dz}Y_M(a,z).
    \]
\end{prp}

\section{The definition of connection twisted modules and their basic properties}\label{s:def and basic properties}

In this section, we introduce the notion of connection twisted modules and investigate their basic properties. The central ingredient of our definition is the formal parallel transport operator $P_\nabla(w,z)$ associated with a $\frg$-connection $\nabla$, which is used to deform the Jacobi identity.

For a $\bbC$-graded vertex algebra $V$, the automorphism group $\ul{\Aut}(V)$ is a functor from the category of commutative $\bbC$-algebras to the category of groups defined by
\[
    \ul{\Aut}(V)(R)\ceq\Aut_R(V_R)=\{g\in\GL_R(V_R)\mid g \text{ is a grade-preserving automorphism on } V_R\}
\]
for any commutative $\bbC$-algebra $R$, where $V_R$ denotes the vertex algebra $V\otimes R$ over $R$. Let $G$ be an affine algebraic group of finite type over $\bbC$. A \emph{vertex algebraic $G$-action} on $V$ is a natural transformation $G\to \ul{\Aut}(V)$. In particular, the group $G(R)$ acts on $V\otimes R$ for a commutative $\bbC$-algebra $R$. This action will be used for $R=\bbC\dpr{z}$ (\cref{s:gauge transformation}) and $R=\wt{\shK}_z$ (\cref{s:regular singular case}). 

To formulate the definition, we first establish the following lemma:
\begin{lem}
    Let $V$ be a $\bbC$-graded vertex algebra with a vertex algebraic $G$-action and let $\nabla$ be a $\frg$-connection $\nabla=d+A(z)dz$ ($A(z)\in\frg\otimes\bbC\dpr{z}$). Then, for any $a,b\in V$, $Y(P_{\nabla}(w,z)a,z)b$ is a well-defined element of $(V\otimes\bbC\dpr{w})\dpr{z}$, where $P_{\nabla}(w,z)$ is a formal power series with coefficients in $U(\frg)\otimes\bbC\dpr{w}$ defined by
    \[
        P_{\nabla}(w,z)=e^{z(\pd_w+A(w))}1=\sum_{s\geq 0}z^sA^{[s]}(w).
    \]
    The first few terms are $A^{[0]}(w)=1$, $A^{[1]}(w)=A(w)$, $A^{[2]}(w)=\frac{A'(w)+A(w)^2}{2}$.
\end{lem}

\begin{proof}
    We explicitly write 
    \[
        Y(P_{\nabla}(w,z)a,z)b=\sum_{s\geq 0}z^sY(A^{[s]}(w)a,z)b=\sum_{n\in\bbZ}\left(\sum_{s\geq 0}(A^{[s]}(w)a)_{(n+s)}b\right)z^{-n-1}.
    \]
     Since the grading of $(A^{[s]}(w)a)_{(n+s)}b$ is $\Delta_a+\Delta_b-n-s-1$, for a fixed $n\in\bbZ$, we have $(A^{[s]}(w)a)_{(n+s)}b=0$ for sufficiently large $s$. Thus, for any $n\in\bbZ$, we have
     \[
        \sum_{s\geq 0}(A^{[s]}(w)a)_{(n+s)}b\in V\otimes\bbC\dpr{w}.
     \]
     By considering the grading again, we have $\sum_{s\geq 0}(A^{[s]}(w)a)_{(n+s)}b=0$ for sufficiently large $n$. Thus, this proves the claim.
\end{proof}

\begin{dfn}\label{dfn:conn tw mod}
    Let $V$ be a $\bbC$-graded vertex algebra with a vertex algebraic $G$-action, and let $\nabla$ be a $\frg$-connection $\nabla=d+A(z)dz$ ($A(z)\in\frg\otimes\bbC\dpr{z}$). A \emph{$\nabla$-twisted $V$-module} $M=(M,Y_M)$ consists of the following data:
    \begin{itemize}
    \item 
    $M$ is a vector superspace,
    \item 
    $Y^\nabla_M\colon V\to\Hom(M,M\dpr{z})$ is a linear map, preserving parities,
    \end{itemize}
    satisfying that
    \begin{enumerate}
    \item 
    $Y^\nabla_M(\one,z)=\id_M$
    \item 
    $\nabla$-twisted Jacobi identity: 
    \begin{equation}\label{eq:nabla tw Jacobi}
        \begin{split}
            &z_0^{-1}\delta\left(\frac{z_1-z_2}{z_0}\right)Y^\nabla_M(a,z_1)Y^\nabla_M(b,z_2)+(-1)^{\abs{a}\abs{b}}(-z_0)^{-1}\delta\left(\frac{z_2-z_1}{-z_0}\right)Y^\nabla_M(b,z_2)Y^\nabla_M(a,z_1)\\
            &=z_2^{-1}\delta\left(\frac{z_1-z_0}{z_2}\right)Y^\nabla_M(Y(P_{\nabla}(z_2,z_0)a,z_0)b,z_2).
        \end{split}
    \end{equation}
    \end{enumerate}
\end{dfn}

Hereafter, let $V$ be a $\bbC$-graded vertex algebra with a vertex algebraic $G$-action and $\nabla$ be a $\frg$-connection $\nabla=d+A(z)dz$ ($A(z)\in\frg\otimes\bbC\dpr{z}$).

As in the usual theory of vertex algebra modules, the twisted Jacobi identity gives rise to a family of identities obtained by taking formal residues. We first derive the following analogue of the Borcherds identity.
\begin{prp}\label{prp:nabla tw Borcherds identity}
    Let $M$ be a $\nabla$-twisted $V$-module. Then, for any $n\in\bbZ$, the following identity, called the \emph{$\nabla$-twisted Borcherds identity}, holds:
    \begin{equation}\label{eq: nabla tw Borcherds identity}
        \begin{split}
            &(z_1-z_2)^nY^\nabla_M(a,z_1)Y^\nabla_M(b,z_2)-(-1)^{\abs{a}\abs{b}}(z_2-z_1)^nY^\nabla_M(b,z_2)Y^\nabla_M(a,z_1)\\
            &=\sum_{m,s\geq 0}Y^\nabla_M((A^{[s]}(z_2)a)_{(m+s+n)}b,z_2)\pd_2^{(m)}z_1^{-1}\delta(z_2/z_1).
        \end{split}
    \end{equation}
\end{prp}

\begin{proof}
    By taking $\oint dz_0\ z_0^n$ of $\nabla$-twisted Jacobi identity \cref{eq:nabla tw Jacobi}, we have 
    \begin{align*}
        &(z_1-z_2)^nY^\nabla_M(a,z_1)Y^\nabla_M(b,z_2)-(-1)^{\abs{a}\abs{b}}(z_2-z_1)^nY^\nabla_M(b,z_2)Y^\nabla_M(a,z_1)\\
        &=\oint dz_0\sum_{m,s\geq 0}z_0^{n+m+s}Y^\nabla_M(Y(A^{[s]}(z_2)a,z_0)b,z_2)\pd_2^{(m)}z_1^{-1}\delta(z_2/z_1)\\
        &=\sum_{m,s\geq 0}Y^\nabla_M((A^{[s]}(z_2)a)_{(n+m+s)}b,z_2)\pd_2^{(m)}z_1^{-1}\delta(z_2/z_1),
    \end{align*}
    where the first equality follows from \cref{eq: skew symmetry} and the definition of $A^{[s]}(z_2)$. 
\end{proof}

The following two corollaries immediately follow from the $\nabla$-twisted Borcherds identity.
\begin{cor}\label{cor: comm formula}
    \begin{equation}\label{eq: comm formula}
        [Y^\nabla_M(a,z_1),Y^\nabla_M(b,z_2)]=\sum_{m,n\geq 0}Y^\nabla_M((A^{[n]}(z_2)a)_{(m+n)}b,z_2)\pd_2^{(m)}z_1^{-1}\delta(z_2/z_1)
    \end{equation}
\end{cor}

\begin{cor}\label{cor: n-th product}
    \begin{equation}\label{eq: n-th product}
        Y^\nabla_M(a,z)_{(n)}Y^\nabla_M(b,z)=\sum_{s\geq 0}Y^\nabla_M((A^{[s]}(z)a)_{(s+n)}b,z).
    \end{equation}
    In particular, we have
    \begin{equation}\label{eq: normally ordered product}
        \nop{Y^\nabla_M(a,z)Y^\nabla_M(b,z)}=Y^\nabla_M(a_{(-1)}b,z)+Y^\nabla_M((A(z)a)_{(0)}b,z)+Y^\nabla_M((A^{[2]}(z)a)_{(1)}b,z)+\cdots.
    \end{equation}
\end{cor}

\begin{prp}
    The axiom of $\nabla$-twisted Jacobi identity (2) is equivalent to the following axioms:
    \begin{enumerate}[label=(2-\arabic*)]
        \item
        weak commutativity (locality): for any $a,b\in V$, there exists $N\geq 0$
        \begin{equation}\label{eq: weak commutativity}
            (z_1-z_2)^NY^\nabla_M(a,z_1)Y^\nabla_M(b,z_2)=(-1)^{\abs{a}\abs{b}}(z_1-z_2)^NY^\nabla_M(b,z_2)Y^\nabla_M(a,z_1),
        \end{equation}
        \item 
        $\nabla$-twisted weak associativity: for any $a,b\in V$, there exists $N\geq 0$
        \begin{equation}\label{eq: nabla tw weak associativity}
            (z_0+z_2)^NY^\nabla_M(a,z_0+z_2)Y^\nabla_M(b,z_2)=(z_0+z_2)^NY^\nabla_M(Y(P_\nabla(z_2,z_0)a,z_0)b,z_2).
        \end{equation}
    \end{enumerate}
\end{prp}

\begin{proof}
    Suppose that $\nabla$-twisted Jacobi identity (2) holds. Then, weak commutativity (2-1) immediately follows from \cref{eq: comm formula}. Hence, we can apply \cref{prp:fundamental property of n-th product} to $Y^\nabla_M(a,z)$ and $Y^\nabla_M(b,w)$ ($a,b\in V$). By \cref{eq: n-th product}, we have 
    \[
        Y^\nabla_M(a,z)Y^\nabla_M(b,w)
        =\sum_{n<N}(z-w)^{-n-1}\sum_{s\geq 0}Y^\nabla_M((A^{[s]}(w)a)_{(n+s)}b,w)
        =Y^\nabla_M(Y(P_\nabla(w,z-w)a,z-w)b,w),
    \]
    where the right hand side is expanded in $\abs{z}>\abs{w}$, which implies weak associativity (2-2).

    Conversely, suppose that weak commutativity (2-1) and weak associativity (2-2) hold. By weak commutativity (2-1), for any $a,b\in V$, $v\in M$, there exists $l\geq0$ such that
    \[
        (z_1-z_2)^lY^\nabla_M(a,z_1)Y^\nabla_M(b,z_2)v\in M\dpr{z_1}\dpr{z_2}\cap M\dpr{z_2}\dpr{z_1}=M\dbr{z_1,z_2}[z_1^{-1},z_2^{-1}].
    \]
    Thus, we can write
    \begin{equation}\label{eq: def of f}
        Y^\nabla_M(a,z_1)Y^\nabla_M(b,z_2)v=\frac{h(z_1,z_2)}{(z_1-z_2)^lz_1^mz_2^n}
    \end{equation}
    for some formal power series $h(z_1,z_2)$ with coefficients in $M$, where the right hand side is expanded in $\abs{z_1}>\abs{z_2}$. We set 
    \[
        f(z_0,z_1,z_2)\ceq \frac{h(z_1,z_2)}{z_0^lz_1^mz_2^n}\in M\dbr{z_0,z_1,z_2}[z_0^{-1},z_1^{-1},z_2^{-1}].
    \]
    Then the equality \cref{eq: def of f}, the weak commutativity (2-1), and the $\nabla$-twisted weak associativity (2-2) mean
    \begin{align*}
        Y^\nabla_M(a,z_1)Y^\nabla_M(b,z_2)v&=f(z_1-z_2,z_1,z_2)\qquad (\abs{z_1}>\abs{z_2}),\\
        (-1)^{\abs{a}\abs{b}}Y^\nabla_M(b,z_2)Y^\nabla_M(a,z_1)v&=f(-z_2+z_1,z_1,z_2) \qquad (\abs{z_2}>\abs{z_1}),\\
        Y^\nabla_M(Y(P_\nabla(z_2,z_0)a,z_0)b,z_2)v&=f(z_0,z_2+z_0,z_2)\qquad (\abs{z_2}>\abs{z_0}),
    \end{align*}
    respectively. Multiplying the identity \cref{eq: Jacobi identity} to $f(z_0,z_1,z_2)$, we obtain $\nabla$-twisted Jacobi identity \cref{eq:nabla tw Jacobi} by \cref{prp:basic property of formal delta funtion}. 
\end{proof}

To simplify \cref{cor: n-th product}, we introduce $\nabla$-deformed $(n)$-th product as follows:
\begin{dfn}\label{dfn:deformed n-th product of tw mods}
\[
    \begin{split}
    Y^\nabla_M(a,w)_{(n,\nabla)}Y^\nabla_M(b,w)
    &=\oint dz \bigl((z-w)^{n+s}Y^\nabla_M(A^{[s]}_\mathrm{inv}(w))Y^\nabla_M(b,w)-(-w+z)^nY^\nabla_M(b,w)Y^\nabla_M(A^{[s]}_\mathrm{inv}(w)a)\bigr)\\
    &=\sum_{s\geq 0}Y^\nabla_M(A^{[s]}_\mathrm{inv}(x)a,w)_{(n+s)}Y^\nabla_M(b,w)|_{x=w},
    \end{split}
\]
where $A^{[s]}_\mathrm{inv}(w)$ is a formal Laurent series with coefficients in $U(\frg)$ defined to be $P_\nabla(w-z,-z)=\sum_{s\geq 0}z^sA^{[s]}_\mathrm{inv}(w)$. The lower terms are $A^{[0]}_\mathrm{inv}(w)=1$, $A^{[1]}_\mathrm{inv}(w)=-A(w)$, $A^{[2]}_\mathrm{inv}(w)=\frac{A(w)^2-A'(w)}{2}$.
\end{dfn}

\begin{lem}\label{lem:inverse transformation}
    The following holds in $U(\frg)\dpr{w}\dbr{z}$
    \[
        P_\nabla(w-z,-z)P_\nabla(w,z)=P_\nabla(w,z)P_\nabla(w-z,-z)=1,
    \]
    which is equivalent to
    \[
        \sum_{s\geq 0}A^{[s]}_\mathrm{inv}(w)A^{[u-s]}(w)=\sum_{s\geq 0}A^{[u-s]}(w)A^{[s]}_\mathrm{inv}(w)=\delta_{u,0}
    \]
    for $u\geq 0$.
\end{lem}

\begin{proof}
    By definition, we have
    \[
        \begin{split}
        P_\nabla(w-z,-z)P_\nabla(w,z)
        &=(e^{-z(A(w-z)+\pd_w)}1)e^{z(A(w)+\pd_w)}1
        =e^{-z(A(w-z)+\pd_w)}e^{-z\pd_w}e^{z(A(w)+\pd_w)}1\\
        &=e^{-z(A(w-z)+\pd_w)}e^{z(A(w-z)+\pd_w)}1=1.
        \end{split}
    \]
    The other is similar.
\end{proof}

Then, \cref{cor: n-th product} is rewritten as:
\begin{cor}\label{cor:deformed n-th product identity}
    For a $\nabla$-twisted $V$-module $M$, the following identity holds:
    \[
        Y^\nabla_M(a,z)_{(n,\nabla)}Y^\nabla_M(b,z)=Y^\nabla(a_{(n)}b,z).
    \]
\end{cor}

\begin{proof}
    By \cref{cor: n-th product}, we have
    \[
        \begin{split}
        Y^\nabla_M(a,z)_{(n,\nabla)}Y^\nabla_M(b,z)
        &=\sum_{s\geq 0}Y^\nabla_M(A^{[s]}_\mathrm{inv}(x)a,z)_{(n+s)}Y^\nabla_M(b,z)|_{x=z}\\
        &=\sum_{s,t\geq 0}Y^\nabla_M((A^{[t]}(z)A^{[s]}_\mathrm{inv}(z)a)_{(n+s+t)}b,z).
        \end{split}
    \]
    Then, by \cref{lem:inverse transformation}, we have
    \[
        \sum_{s,t\geq 0}(A^{[t]}(z)A^{[s]}_\mathrm{inv}(z)a)_{(n+s+t)}b=a_{(n)}b.
    \]
    Thus, we obtain the claim.
\end{proof}

The translation property for ordinary vertex algebra modules is deformed by the connection.
\begin{cor}
    For a $\nabla$-twisted $V$-module $M$, the following holds:
    \[
        \pd_z Y_M(a,z)=Y_M((A(z)+T)a,z)
    \]
\end{cor}

\begin{proof}
    Applying the $\nabla$-twisted weak associativity \cref{eq: nabla tw weak associativity} with $b=\one$, we have
    \[
        (z_0+z_2)^NY^\nabla_M(a,z_0+z_2)v=(z_0+z_2)^NY^\nabla_M(Y(P_\nabla(z_2,z_0)a,z_0)\one,z_2)v
    \]
    for any $v\in M$. By taking $N$ large enough, we may assume that $(z_0+z_2)^NY^\nabla_M(a,z_0+z_2)v\in M\dbr{z_0,z_2}$. Thus, we have
    \[
        (z_0+z_2)^Ne^{z_0\pd_2}Y^\nabla_M(a,z_2)v=(z_0+z_2)^NY^\nabla_M(Y(P_\nabla(z_2,z_0)a,z_0)\one,z_2)v.
    \]
    Since both $e^{z_0\pd_2}Y^\nabla_M(a,z_2)v$ and $Y^\nabla_M(Y(P_\nabla(z_2,z_0)a,z_0)\one,z_2)v$ are contained in $M\dpr{z_2}\dbr{z_0}$, we have
    \[
        e^{z_0\pd_2}Y^\nabla_M(a,z_2)v=Y^\nabla_M(Y(P_\nabla(z_2,z_0)a,z_0)\one,z_2)v=\sum_{n\geq 0}z_0^nY^\nabla_M\left(\sum_{0\leq k\leq n}T^{(k)}A^{[n-k]}(z_2)a,z_2\right)v.
    \]
    Comparing the coefficients of $z_0$, we obtain the desired formula.
\end{proof}

Since $P_\nabla(w,z)=1$ on $V^{G}$, we obtain the following proposition:  
\begin{prp}\label{prp:rest on G-invariant}
    Let $M$ be a $\nabla$-twisted $V$-module. The restriction of $M$ to $V^{G}$ is a (untwisted) $V^G$-module. In particular, if $V$ has a $G$-invariant conformal element $\omega\in V$, then $L^\nabla_M(z)=Y^\nabla_M(\omega,z)$ yields an action of Virasoro algebra on $M$. 
\end{prp}

\section{The gauge transformations}\label{s:gauge transformation}
We now investigate the behavior of connection twisted modules under gauge transformations. Since gauge-equivalent connections describe the same underlying connection in different trivializations, it is natural to expect that their associated module categories are equivalent. The following theorem confirms this expectation.

For a $\frg$-connection $\nabla=d+A(z)dz$ and $g=g(z)\in G(\bbC\dpr{z})$, we define a $\frg$-connection $\nabla^g$ by
\[
    \nabla^g\ceq g(z)^{-1}\nabla g(z)=d+g(z)^{-1}dg(z)+g(z)^{-1}A(z)g(z)dz.
\]
This transformation $\nabla\mto\nabla^g$ is called a \emph{gauge transformation} by $g$.

\begin{thm}\label{thm:gauge transf}
    Let $g=g(z)\in G(\bbC\dpr{z})$. For a $\nabla$-twisted module $M=(M,Y^\nabla_M)$, we define a field $(Y^\nabla_M)^g(a,z)$ by $(Y^\nabla_M)^g(a,z)\ceq Y^\nabla_M(g(z)a,z)$ for $a\in V$. Then, $M^g=(M,(Y^\nabla_M)^g)$ is a $\nabla^g$-twisted $V$-module.
\end{thm}

\begin{proof}
    To begin with, we establish the well-definedness of $Y^\nabla_M(g(z)a,z)$ and verify weak commutativity. Let $g(z)a=\sum_{i=1}^ka_i\otimes g^i_a(z)$, $a_i\in V$, $g_a^i(z)\in\bbC\dpr{z}$. Then, $(Y^\nabla_M)^g(a,z)=Y^\nabla_M(g(z)a,z)$ is written as $\sum_{i=1}^kY^\nabla_M(a_i,z)g^i_a(z)$. Thus, we find that $Y^\nabla_M(g(z)a,z)$ is well-defined element of $\Hom(M,M\dpr{z})$ and 
    \[
        (z-w)^N[(Y^\nabla_M)g(a,z),(Y^\nabla_M)^g(b,w)]=\sum_{i,j}(z-w)^N[Y^\nabla_M(a_i,z),Y^\nabla_M(b_j,w)]g^i_a(z)g^j_b(w)=0
    \]
    for sufficiently large $N\geq 0$. Since $g(z)\one=\one$, we have $(Y^\nabla_M)^g(\one,z)=\id_M$. It remains to show $\nabla^g$-twisted weak associativity. For any $a,b\in V$, there exists $N\geq 0$ such that
    \begin{align*}
        (z_0+z_2)^N(Y^\nabla_M)^g(a,z_0+z_2)(Y^\nabla_M)^g(b,z_2)
        &=(z_0+z_2)^NY^\nabla_M(g(z_0+z_2)a,z_0+z_2)Y^\nabla_M(g(z_2)b,z_2)\\
        &=(z_0+z_2)^NY^\nabla_M(Y(P_\nabla(z_2,z_0)g(z_2+z_0)a,z_0)g(z_2)b,z_2).
    \end{align*}
    Here, we have
    \[
        P_\nabla(z_2,z_0)g(z_2+z_0)=(e^{z_0(\pd_2+A(z_2))}1)e^{z_0\pd_2}g(z_2)=e^{z_0(\pd_2+A(z_2))}g(z_2)=g(z_2)P_{\nabla^g}(z_2,z_0).
    \]
    Thus, we obtain the $\nabla^g$-twisted weak associativity:
    \begin{align*}
        (z_0+z_2)^N(Y^\nabla_M)^g(a,z_0+z_2)(Y^\nabla_M)^g(b,z_2)
        &=(z_0+z_2)^NY^{\nabla}_M(Y(g(z_2)P_{\nabla^g}(z_2,z_0)a,z_0)g(z_2)b,z_2)\\
        &=(z_0+z_2)^NY^{\nabla}_M(g(z_2)Y(P_{\nabla^g}(z_2,z_0)a,z_0)b,z_2)\\
        &=(z_0+z_2)^N(Y^\nabla_M)^g(Y(P_{\nabla^g}(z_2,z_0)a,z_0)b,z_2).
    \end{align*}
\end{proof}

We denote by $\cMod^\nabla(V)$ the category of $\nabla$-twisted $V$-modules, whose objects are $\nabla$-twisted $V$-modules and whose morphisms are linear maps commuting with the action of $V$. This construction defines a functor $(-)^g$ from $\cMod^\nabla(V)$ to $\cMod^{\nabla^g}(V)$. By construction, we have
\[
    (-)^{gh}\simeq (-)^{h}\circ (-)^{g},\quad (-)^{e}\simeq\id
\]
for $g,h\in G(\bbC\dpr{z})$ and the identity element $e\in G(\bbC\dpr{z})$. Thus, we have the following:
\begin{cor}\label{cor:main A}
    For $g\in G(\bbC\dpr{z})$, the functor $(-)^g\colon\cMod^\nabla(V)\to \cMod^{\nabla^g}(V)$ is an equivalence of categories with a quasi-inverse $(-)^{g^{-1}}$.
\end{cor}

\section{The regular singular case}\label{s:regular singular case}

In this section, we study connection twisted modules associated with regular singular connections. Our main result is that the corresponding module category is equivalent to the category of logarithmic twisted modules associated with the monodromy automorphism. We begin by recalling the definition of logarithmic twisted modules \cite{Hua10,Bak16} in a form convenient for our purposes.

\subsection{Logarithmic $g$-twisted modules}

    We define a commutative algebra $\wt{\shK}_z$ by
    \[
        \wt{\shK}_z=\bigoplus_{\alpha\in\bbC/\bbZ,m\geq 0}\bbC\dpr{z}z^\alpha(\log z)^m
    \]
    and $\pd_z$ is the derivation on $\wt{\shK}_z$ defined by
    \[
        \pd_zz^{\alpha}=\alpha z^{\alpha-1},\quad \pd_z\log z=\frac{1}{z}.
    \]

\begin{dfn}
    For any vector superspace $M$, we set 
    \[
        \LFie(M)\ceq \bigoplus_{\alpha\in\bbC/\bbZ,m\geq 0}\Hom(M,M\dpr{z})z^{\alpha}(\log z)^m.
    \]
    We call an element of $\LFie(M)$ a \emph{logarithmic field} on $M$. If logarithmic fields $a(z)$ and $b(w)$ satisfy that $(z-w)^N[a(z),b(w)]=0$ for $N\gg 0$, then $a(z)$ and $b(w)$ are \emph{local}. Since $(z-w)^Na(z)b(w)$ is regular on $z=w$ for sufficiently large $N\geq 0$, we can define a \emph{$(n)$-th product} $a(w)_{(n)}b(w)\in\LFie(M)$ by
    \[
        a(w)_{(n)}b(w)\ceq \pd_z^{(N-n-1)}(z-w)^Na(z)b(w)|_{z=w}
    \]
    for $n<N$ and $a(w)_{(n)}b(w)=0$ for $n\geq N$. We easily see that the $(n)$-th product is independent of $N$.
\end{dfn}

\begin{rmk}
    If local logarithmic fields $a(z)$ and $b(w)$ are fields, then the above $(n)$-th product $a(w)_{(n)}b(w)$ coincides with that in \cref{dfn: local field and n-th product}.
\end{rmk}

\begin{dfn}
    Let $V$ be a vertex algebra and let $g$ be an automorphism on $V$, whose action on $V$ is locally finite.  A \emph{$g$-twisted $V$-module} $M=(M,Y^g_M)$ consists of the following data:
    \begin{itemize}
        \item 
        $M$ is a vector superspace,
        \item 
        $Y^g_M\colon V\to \LFie(M)$ is a linear map, preserving parities,
    \end{itemize}
    satisfying that
    \begin{enumerate}
        \item 
        $e^{2\pi iz\pd_z}Y^g_M(a,z)=Y^g_M(ga,z)$
        \item 
        $Y^g_M(\one,z)=\id_M$,
        \item
        for any $a,b\in V$, $Y^g_M(a,z)$, $Y^g_M(b,w)$ are local,
        \item 
        $Y^g_M(a,z)_{(n)}Y^g_M(b,z)=Y^g_M(a_{(n)}b,z)$ for any $a,b\in V$, $n\in\bbZ$.
    \end{enumerate}
\end{dfn}

\begin{rmk}
    This definition may appear to be stronger than that of \cite{Bak16}, since the space of logarithmic fields $\LFie(M)$ is smaller than the corresponding space considered there. However, since $g$ acts locally finitely on $V$, the two definitions are equivalent by \cite[Lemma~4.8]{Bak16}.  
\end{rmk}

\subsection{Regular singular connection twisted modules}
In this subsection, let $\nabla=d+A(z)dz$ be a regular singular $\frg$-connection, $F=F(z)$ be an element of $G(\wt{\shK}_z)$ satisfying $F(z)^{-1}\nabla F(z)=d$ (such an element exists, see in \cref{thm:existence of G-valued solutions}), and $g\in G(\bbC)$ be the monodromy of $F(z)$, i.\ e.\ $e^{2\pi i z\pd_z}F(z)=F(z)g$.

The construction in this subsection is essentially an extension of the functor $(-)^g$ ($g\in G(\bbC\dpr{z})$) introduced in \cref{s:gauge transformation} to a fundamental solution $F(z)\in G(\wt{\shK}_z)$. However, unlike an ordinary gauge transformation in $G(\mathbb C((z)))$, the fundamental solution $F(z)$ may be multivalued. Consequently, the resulting fields need not be ordinary fields, and their multivaluedness must be treated carefully. As we shall see, the monodromy of $F(z)$ gives rise precisely to the twisting of a logarithmic twisted module.

    \begin{thm}\label{thm:conn tw mod to log tw mod}
        For a $\nabla$-twisted module $M=(M,Y^\nabla_M)$, we define a logarithmic field $(Y^\nabla_M)^F(a,z)\in\LFie(M)$ by $(Y^\nabla_M)^F(a,z)\ceq Y^\nabla_M(F(z)a,z)$ for $a\in V$. Then, $M^F=(M,(Y^\nabla_M)^F)$ is a $g$-twisted $V$-module.
    \end{thm}

    \begin{proof}
    Let $F(z)a=\sum_{i=1}^ka_i\otimes F_a^i(z)$, $a_i\in V$, $F^i_a(z)\in\wt{\shK}_z$. Then, $(Y^\nabla_M)^F(a,z)=Y^\nabla_M(F(z)a,z)$ is written as $\sum_{i=1}^kY^\nabla_M(a_i,z)F_a^i(z)$. Thus, we find that $(Y^\nabla_M)^F(a,z)$ belongs to $\LFie(M)$, and $(Y^\nabla_M)^F(a,z)$ ($a\in V$) are local. By  $e^{2\pi iz\pd_z}F(z)=F(z)g$, we have $\sum_{i=1}^ka_i\otimes e^{2\pi iz\pd_z}F^i_a(z)=\sum_{i=1}^k(ga)_i\otimes F^i_{ga}(z)$. Thus, it follows that
    \begin{align*}
        e^{2\pi iz\pd_z}(Y^\nabla_M)^F(a,z)
        &=\sum_{i=1}^kY^\nabla_M(a_i,z)e^{2\pi iz\pd_z}F^i_a(z)
        =\sum_{i=1}^kY^\nabla_M((ga)_i,z)F^i_{ga}(z)\\
        &=\sum_{i=1}^kY^\nabla_M(F(z)ga,z)
        =(Y^\nabla_M)^F(ga,z).
    \end{align*}

    It remains to verify the $(n)$-th product identity. By $\nabla$-twisted weak associativity, we have
    \begin{align*}
        (Y^\nabla_M)^F(a,z)(Y^\nabla_M)^F(b,w)
        &=Y^\nabla_M(F(z)a,z)Y^\nabla_M(F(w)b,w)\\
        &=Y^\nabla_M(Y(P_\nabla(w,z-w)F(z)a,z-w)F(w)b,w)v\\
        &=\sum_{n<N}\sum_{s\geq 0}(z-w)^{-n-1}Y^\nabla_M((A^{[s]}(w)F(z)a)_{(n+s)}F(w)b,w),
    \end{align*} 
    in which $(z-w)^{-n-1}$ is expanded in $\abs{z}>\abs{w}$. By multiplying $(z-w)^N$ and taking $\pd_z^{(N-m-1)}$, we have
    \begin{align*}
        &\pd_z^{(N-m-1)}(z-w)^N(Y^\nabla_M)^F(a,z)(Y^\nabla_M)^F(b,w)\\
        &=\sum_{n<N}\sum_{s\geq 0}\sum_{k\geq 0}\binom{N-n-1}{N-m-k-1}(z-w)^{-n+m+k}Y^\nabla_M((A^{[s]}(w)\pd_z^{(k)}F(z)a)_{(n+s)}F(w)b,w).
    \end{align*}
    Substituting $w$ for $z$, we obtain
    \begin{align*}
        (Y^\nabla_M)^F(a,w)_{(n)}(Y^\nabla_M)^F(b,w)
        &=\sum_{s,k\geq 0}Y^\nabla_M((A^{[s]}(w)\pd_w^{(k)}F(w)a)_{(n+k+s)}F(w)b,w)\\
        &=\sum_{s\geq 0}Y^\nabla_M\left(\left(\sum_{k=0}^sA^{[s-k]}(w)\pd_w^{(k)}F(w)a\right)_{(n+s)}F(w)b,w\right).
     \end{align*}
     Since
     \[
        P_\nabla(w,z)F(w+z)=(e^{z(\pd_w+A(w))}1)e^{z\pd_w}F(w)=e^{z(\pd_w+A(w))}F(w)=F(w),
     \]
     which is equivalent to  
     \[
        \sum_{k=0}^sA^{[s-k]}(w)\pd^{(k)}_wF(w)=\delta_{s,0}F(w)
    \]
    for any $s\geq 0$, we have the $(n)$-th product identity:
     \[
        (Y^\nabla_M)^F(a,w)_{(n)}(Y^\nabla_M)^F(b,w)=Y^\nabla_M((F(w)a)_{(n)}F(w)b,w)=Y^\nabla_M(F(w)a_{(n)}b,w)=(Y^\nabla_M)^F(a_{(n)}b,w).
     \]
     \end{proof}

Conversely, we extend the gauge transformation by $F(z)^{-1}$ to logarithmic twisted modules. Although $F(z)^{-1}$ is generally multivalued, its monodromy cancels the twisting of the logarithmic fields. As a result, the fields obtained in this way are ordinary fields and satisfy the $\nabla$-twisted Jacobi identity.
\begin{thm}\label{thm:log tw mod to conn tw mod}
    For a logarithmic $g$-twisted module $M=(M,Y^g_M)$, we define a field $(Y^g_M)^{F^{-1}}(a,z)\in\Hom(M,M\dpr{z})$ by $(Y^g_M)^{F^{-1}}(a,z)\ceq Y^g_M(F(z)^{-1}a,z)$ for $a\in V$. Then, $M^{F^{-1}}=(M,(Y^g_M)^{F^{-1}})$ is a $\nabla$-twisted $V$-module.
\end{thm}

\begin{proof}
    We first verify that $(Y^g_M)^{F^{-1}}(a,z)=Y^g_M(F(z)^{-1}a,z)$ is contained in $\Hom(M,M\dpr{z})$. Let $F(z)^{-1}a=\sum_{i=1}^ka_i\otimes G^i_a(z)$, $a_i\in V$, $G^i_a(z)\in\wt{\shK}_z$. Since $(Y^g_M)^{F^{-1}}(a,z)=Y^g_M(F(z)^{-1}a,z)$ is written as $\sum_{i=1}^kY^g_M(a_i,z)G_a^i(z)$, we have $(Y^g_M)^{F^{-1}}(a,z)\in\LFie(M)$. Note that $\sum_{i=1}^ka_i\otimes e^{2\pi iz\pd_z}G_a^i(z)=\sum_{i=1}^kg^{-1}a_i\otimes G^i_a(z)$ by $e^{2\pi iz\pd_z}F(z)^{-1}=g^{-1}F(z)^{-1}$. By the $g$-equivariance of $Y^g_M$, we have
    \begin{align*}
        e^{2\pi iz\pd_z}(Y^g_M)^{F^{-1}}(a,z)
        &=\sum_{i=1}^k(e^{2\pi iz\pd_z}Y^g_M(a_i,z))e^{2\pi iz\pd_z}G^i_a(z)
        =\sum_{i=1}^kY^g_M(gg^{-1}a_i,z)G^i_a(z)\\
        &=\sum_{i=1}^kY^g_M(a_i,z)G^i_a(z)
        =(Y^g_M)^{F^{-1}}(a,z).
    \end{align*}
    Since $e^{2\pi iz\pd_z}z^{\alpha}(\log z)^m=e^{2\pi i\alpha}z^{\alpha}(\log z+2\pi i)^m$, we have $(Y^g_M)^{F^{-1}}(a,z)\in \Hom(M,M\dpr{z})$.

    Since $(Y^g)^{F^{-1}}(\one,z)=\id_M$ and weak commutativity is immediate, it remains to verify $\nabla$-twisted weak associativity. 
    \begin{align*}
        &\pd_z^{(N-n-1)}(z-w)^N(Y^g_M)^{F^{-1}}(a,z)(Y^g_M)^{F^{-1}}(b,w)\\
        &=\sum_{i,j}\pd_z^{(N-n-1)}(z-w)^NY^g_M(a_i,z)Y^g_M(b_j,w)G^i_a(z)G^j_b(w)\\
        &=\sum_{i,j}\sum_{s\geq 0}\bigl(\pd_z^{(N-n-s-1)}(z-w)^NY^g_M(a_i,z)Y^g_M(b_j,w)\bigr)\pd_z^{(s)}G^i_a(z)G^j_b(w)
    \end{align*}
    Substituting $w$ for $z$, we have
    \begin{align*}
        (Y^g_M)^{F^{-1}}(a,w)_{(n)}(Y^g_M)^{F^{-1}}(b,w)
        &=\sum_{i,j}\sum_{s\geq 0}Y^g_M(a_i,w)_{(n+s)}Y^g_M(b_j,w)(\pd_w^{(s)}G^i_a(w))G^j_b(w)\\
        &=\sum_{s\geq 0}Y^g_M(\pd^{(s)}_wF(w)^{-1}a,w)_{(n+s)}Y^g_M(F(w)^{-1}b,w)
    \end{align*}
    By the $(n)$-th product identity, we have
    \begin{align*}
        (Y^g_M)^{F^{-1}}(a,z)_{(n)}(Y^g_M)^{F^{-1}}(b,w)
        =\sum_{s\geq 0}Y^g_M((\pd_w^{(s)}F(w)^{-1}a)_{(n+s)}F(w)^{-1}b,w)
    \end{align*}
    Here, we have
    \[
        e^{z\pd_w}F(w)^{-1}=F(w)^{-1}e^{z(\pd_w+A(w))}1=F(w)^{-1}P_{\nabla}(w,z),
    \]
    which is equivalent to 
    \[
        \pd_w^{(s)}F(w)^{-1}=F(w)^{-1}A^{[s]}(w)
    \]
    for any $s\geq 0$. Thus, we obtain
    \begin{align*}
        (Y^g_M)^{F^{-1}}(a,z)_{(n)}(Y^g_M)^{F^{-1}}(b,w)
       &=\sum_{s\geq 0}Y^g_M((F(w)^{-1}A^{[s]}(w)a)_{(n+s)}F(w)^{-1}b,w)\\
       &=\sum_{s\geq 0}Y^g_M(F(w)^{-1}(A^{[s]}(w)a)_{(n+s)}b,w)\\
       &=\sum_{s\geq 0}(Y^g_M)^{F^{-1}}((A^{[s]}(w)a)_{(n+s)}b,w).
    \end{align*}
    Since we have already seen that $(Y^g_M)^{F^{-1}}(a,z)\in\Hom(M,M\dpr{z})$, the above identity is equivalent to $\nabla$-twisted weak associativity:
    \[
        (z_0+z_2)^N(Y^g_M)^{F^{-1}}(a,z_0+z_2)(Y^g_M)^{F^{-1}}(b,z_2)=(z_0+z_2)^N(Y^g_M)^{F^{-1}}(Y(P_\nabla(z_2,z_0)a,z_0)b,z_2).
    \]
\end{proof}

We denote by $\cMod^g(V)$ the category of $g$-twisted $V$-modules, whose objects are $g$-twisted $V$-modules and whose morphisms are linear maps commuting with the action of $V$. Then, by combining the above theorems, we obtain the following corollary:
\begin{cor}\label{cor:main B}
    Let $V$ be a $\bbC$-graded vertex algebra with vertex algebraic $G$-action. Let $\nabla=d+A(z)dz$ be a regular singular $\frg$-connection, $F(z)\in G(\wt{\shK}_z)$ be an element such that $F(z)^{-1}\nabla F(z)=d$, and $g\in G(\bbC)$ be a monodromy of $F(z)$. Then, the functor
    \[
        (-)^F\colon \cMod^\nabla(V)\lto \cMod^g(V)
    \]
    is an equivalence of categories with its quasi-inverse $(-)^{F^{-1}}$.
\end{cor}
This equivalence shows that, in the regular singular case, the category of connection twisted modules is completely determined, up to equivalence, by the monodromy automorphism. In particular, the theory of connection twisted modules recovers the usual theory of logarithmic twisted modules for regular singular connections, which is expected in \cite{FL25,FL26}

\section{Examples}\label{s:example}
In this section, we apply the theory developed above to the symplectic fermion and affine vertex algebras.

\subsection{The $\nabla$-twisted modules of symplectic fermion}
In \cite{FL25}, the authors studied a deformation of its mode algebra associated with an $\frsl_2$-valued connection and constructed a Virasoro action on its modules. We show that both structures arise naturally from the general theory of connection twisted modules.

We consider the $2$-dimensional vector space $\bbC^2=\bbC e_+\oplus\bbC e_-$ with the standard symplectic form $(,)$:
\[
    \begin{bmatrix}
        (e_+,e_+) & (e_+,e_-)\\
        (e_-,e_+) & (e_-,e_-)
    \end{bmatrix}
    =
    \begin{bmatrix}
        0 & 1 \\
        -1 & 0
    \end{bmatrix}
    .
\]
Using these data, we define the symplectic fermion:
\begin{dfn}
    The symplectic fermion $\mathcal{SF}$ is a vertex algebra strongly generated by $\psi_+(z),\psi_-(w)$, both of parity $1$ and grading $1$, with the OPE given by
    \[
        \psi_a(z)\psi_b(w)\sim \frac{(e_a,e_b)}{(z-w)^2}
    \]
for $a,b\in\{+,-\}$. The symplectic fermion $\mathcal{SF}$ has a conformal vector $\omega=\psi_{-,-1}\psi{+,-1}\one$ at central charge $-2$.
\end{dfn}

In what follows, we write $\psi_+,\psi_-\in\mathcal{SF}$ for $\psi_{+,-1}\one,\psi_{-,-1}\one$, respectively.
$\mathcal{SF}$ has a vertex algebraic $\SL_2$-action defined by
\[
    g\psi_+\ceq a\psi_++c\psi_-,
    \quad
    g\psi_-\ceq b\psi_++d\psi_-,
    \quad
    g=
    \begin{bmatrix}
        a & b\\
        c & d
    \end{bmatrix}
    \in\SL_2.
\]

Let $\nabla=d+A(z)dz$ ($A(z)\in\frsl_2\otimes\bbC\dpr{z}$) and $M$ be a $\nabla$-twisted $\mathcal{SF}$-module. Setting $\psi^\nabla_a(z)\ceq Y^\nabla_M(\psi_a,z)$, we have the commutator formula by \cref{cor: comm formula}:
\begin{equation}\label{eq:deformed mode algeba of symplectic algebra}
        [\psi_a(z),\psi_b(w)]
        =(e_a,e_b)\pd_ww^{-1}\delta(z/w)+C_{ab}(w)w^{-1}\delta(z/w),
\end{equation}
where
\begin{equation}\label{eq:def of C}
    C_{ab}(w)=\sum_{k}C_{ab,k}w^{-k-1}\ceq(A(w)e_a,e_b).
\end{equation}

Thus, we have the following claim:
\begin{prp}\label{prp:deforming mode algebra of symp ferm}
    The operators $\psi_{a,m}^\nabla$ ($a\in\{+,-\}$, $m\in\bbZ$), defined to be $\psi^\nabla_a(z)=\sum_{m\in\bbZ}\psi^\nabla_{a,m}z^{-m-1}$, forms the deformed mode algebra introduced in \cite[Definition~2.1]{FL25}.  
\end{prp}

We now examine the conformal structure of connection twisted modules over symplectic fermion. Since the conformal vector $\omega=\psi_{-,-1}\psi_+$ is invariant under the $SL_2$-action, \cref{prp:rest on G-invariant} provides a Virasoro action at central charge $-2$ on every $\nabla$-twisted $\mathcal{SF}$-module. By \cref{cor: n-th product}, we can write explicitly:
\begin{align*}
    L^\nabla(z)
    &=\nop{\psi^\nabla_-(z)\psi^\nabla_+(z)}-\frac{1}{2}((A'(z)+A(z)^2)e_-,e_+)\\
    &=\nop{\psi^\nabla_-(z)\psi^\nabla_+(z)}-\frac{1}{2}(C_{-+})'(z)-\frac{1}{2}(C_{++}(z)C_{--}(z)-C_{+-}(z)C_{-+}(z)),
\end{align*}
that is
\begin{equation}\label{eq:conformal str}
    L^\nabla_n=\sum_{i\in\bbZ}\nop{\psi^\nabla_{-,i}\psi^\nabla_{+,n-i}}+\frac{n+1}{2}C_{-+,n}-\frac{1}{2}\sum_{i}(C_{++,i}C_{--,n-i}-C_{+-,i}C_{-+,n-i}).
\end{equation}
Here, $\nop{\psi^\nabla_{a,i}\psi^\nabla_{b,j}}$ is defined by
\[
    \nop{\psi^\nabla_{a,i}\psi^\nabla_{b,j}}=\left\{
    \begin{array}{cl}
        \psi^\nabla_{a,i}\psi^\nabla_{b,j} & (i<0)\\
        -\psi^\nabla_{b,j}\psi^\nabla_{a,i} & (i\geq 0).
    \end{array}
    \right.
\]

On the other hand, in \cite[\S 4]{FL25}, the authors construct a Virasoro action on a module of the deformed mode algebra of symplectic fermion. Here, we review this construction. For $\psi^\nabla_{a,i}$, we define $:\!\psi^\nabla_{a,i}\psi^\nabla_{b,j}\!:$ by
\[  
    :\!\psi^\nabla_{a,i}\psi^\nabla_{b,j}\!:=\left\{
    \begin{array}{cl}
        \psi^\nabla_{a,i}\psi^\nabla_{b,j} & (i<j)\\
        -\psi^\nabla_{b,j}\psi^\nabla_{a,i} & (i>j)\\
        \frac{1}{2}(\psi^\nabla_{a,i}\psi^\nabla_{b,j}-\psi^\nabla_{b,j}\psi^\nabla_{a,i}) & (i=j).
    \end{array}
    \right.
\]
Then, by \cite[Lemma~4.3]{FL25}, the following operators\footnote{There is a sign error in \cite[Lemma~4.3]{FL25}. The corrected formula is given by \cref{eq:conformal str in FL}.} define an action of Virasoro algebra at central charge $c=-2$:
\begin{equation}\label{eq:conformal str in FL}
    \wt{L}^\nabla_n=\sum_{i\in\bbZ}:\!\psi^\nabla_{-,i}\psi^\nabla_{+,n-i}\!:+c_n,\quad c_n=\frac{1}{2}\sum_{i\in\bbZ}(C_{-+,i}C_{+-,n-i}-C_{++,i}C_{--,n-i}).
\end{equation}

At first sight, the two expressions for the Virasoro operators appear to differ because they use different normal ordering conventions. The difference between these conventions is precisely accounted for by the correction term involving the connection. We make this comparison explicit in the following proposition:
\begin{prp}
    For any $\nabla$-twisted $\mathcal{SF}$-module $M$, the operator $L^\nabla_n$ defined by \cref{eq:conformal str} coincides with $\wt{L}^\nabla_n$ for any $n\in\bbZ$. 
\end{prp}

\begin{proof}
    It suffices to show that 
    \[
        \sum_{i\in\bbZ}:\!\psi^\nabla_{a,i}\psi^\nabla_{b,n-i}\!:=\sum_{i\in\bbZ}\nop{\psi^\nabla_{a,i}\psi^\nabla_{b,n-i}}+\frac{n+1}{2}C_{ab,n}.
    \]    
    We only check the case of $n\geq 0$ here, as the other cases can be checked in the same way. If $n$ is odd, then
    \begin{align*}
        \sum_{i\in\bbZ}:\!\psi^\nabla_{a,i}\psi^\nabla_{b,n-i}\!:
        &=\sum_{i\in\bbZ}\nop{\psi^\nabla_{a,i}\psi^\nabla_{b,n-i}}+\sum_{0\leq i\leq (n-1)/2}[\psi^\nabla_{a,i},\psi^\nabla_{b,n-i}]\\
        &=\sum_{i\in\bbZ}\nop{\psi^\nabla_{a,i}\psi^\nabla_{b,n-i}}+\frac{n+1}{2}C_{ab,n}.
    \end{align*}
    If $n$ is even, then
    \begin{align*}
        \sum_{i\in\bbZ}:\!\psi^\nabla_{a,i}\psi^\nabla_{b,n-i}\!:
        &=\sum_{i\in\bbZ}\nop{\psi^\nabla_{a,i}\psi^\nabla_{b,n-i}}+\sum_{0\leq i\leq n/2-1}[\psi^\nabla_{a,i},\psi^\nabla_{b,n-i}]+\frac{1}{2}[\psi^\nabla_{a,n/2},\psi^\nabla_{b,n/2}]\\
        &=\sum_{i\in\bbZ}\nop{\psi^\nabla_{a,i}\psi^\nabla_{b,n-i}}+\frac{n+1}{2}C_{ab,n},
    \end{align*}
    where $\sum_{0\leq i\leq n/2-1}=0$ if $n=0$.
\end{proof}

\subsection{$\nabla$-twisted modules of affine vertex algebras}\label{ss:affine va}
\begin{dfn}
    Let $\frg$ be a simple Lie algebra, $(,)$ be the standard bilinear form, and $k$ be a complex number. The affine vertex algebra $V^k(\frg)$ is a vertex algebra strongly generated by $x(z)$ ($x\in\frg$), parity $0$ and grading $1$, with OPE given by
    \[
        x(z)y(w)\sim \frac{[x,y](w)}{z-w}+\frac{k(x,y)}{(z-w)^2}
    \]
    for $x,y\in\frg$. If $k$ is not the minus dual Coxeter number $-h^\vee$, then the affine vertex algebra $V^k(\frg)$ has a conformal vector 
    \[
        \omega_\mathrm{sug}=\frac{x_{i,-1}x^i_{-1}\one}{2(k+h^\vee)}
    \]
    at central charge $k\dim\frg/(k+h^\vee)$, where $x_i$ is a basis of $\frg$, $x^i\in\frg$ be its dual basis with respect to $(,)$.   
\end{dfn}

Since $\frg$ acts on $V^k(\frg)$ as $xa\ceq x_0a$ for $x\in\frg$, $a\in V^k(\frg)$, we consider $\nabla$-twisted $V^k(\frg)$-modules for a $\frg$-valued connection $\nabla=d+A(z)dz$. 

In what follows, let $x$ be an element $x_{-1}\one$ ($x\in\frg$) of $V^k(\frg)$ and $x^\nabla(z)=Y^\nabla_M(x,z)$ for $\nabla$-twisted $V^k(\frg)$-module $M$. Then, by \cref{cor: comm formula}, we have
\[
    [x^\nabla(z),y^\nabla(w)]=\bigl([x,y]^\nabla(w)+k(A(w),[x,y])\bigr)w^{-1}\delta(z/w)+k(x,y)\pd_ww^{-1}\delta(z?w)
\]
for $x,y\in\frg$. Thus, we have the following:
\begin{prp}\label{prp:deform mode alg aff va}
    We define a Lie algebra $\wh{\frg}_\nabla$ by $\wh{\frg}_\nabla=\frg\dpr{t}\oplus\bbC K$ with the relation given by
    \[
        [x^\nabla_m,y^\nabla_n]=[x,y]^\nabla_{m+n}+(A_{m+n},[x,y])K+m\delta_{m+n,0}(x,y)K
    \]
    for $x^\nabla_m=xt^m$ and $K$ is central. Then, the mode $x^\nabla_m$ of $x^\nabla(z)$ form the Lie algebra $\wh{\frg}_\nabla$ at level $k$.
\end{prp}

Since the map $(x_m,y_n)\mto (A_{m+n},[x,y])$ is a $2$-coboundary of $\frg\dpr{t}$, $\wh{\frg}_\nabla$ is isomorphic to the affine Lie algebra $\wh{\frg}$ as Lie algebras, that is
\begin{prp}\label{prp:isom aff Lie alg}
    The map $\wh{\frg}\to\wh{\frg}_\nabla$ given by
    \[
        x_n\lmto x^\nabla_n+(A_n,x)K,
        \quad
        K\lmto K,
    \]
    for $x\in\frg$, $n\in\bbZ$, is an isomorphism as Lie algebras.
\end{prp}

Thus, if we prove that the category $\cMod^\nabla(V^k(\frg))$ of $\nabla$-twisted $V^k(\frg)$-modules is equivalent to the category of continuous $\wh{\frg}_\nabla$-modules, then we obtain a categorical equivalence $\cMod(V^k(\frg))\simeq\cMod^\nabla(V^k(\frg))$, which is expected in \cite{FL25,FL26}.

Finally, we determine the Virasoro action $L^\nabla_\mathrm{sug}(z)=Y^\nabla_M(\omega_\mathrm{sug},z)$ on a $\nabla$-twisted $V^k(\frg)$-module $M$. By \cref{cor: n-th product}, we have
\[
    \begin{split}
    L^\nabla_\mathrm{sug}(z)
    &=\frac{1}{2(k+h^\vee)}\sum_i
    \left(
    \nop{x^\nabla_i(z)x^{\nabla,i}(z)}-[[A(z),x_i],x^i]^\nabla(z)-\frac{k}{2}([A'(z),x_i]+[A(z),[A(z),x_i]],x^i)
    \right)\\
    &=\frac{1}{2(k+h^\vee)}\sum_i
    \bigl(
    \nop{x^\nabla_i(z)x^{\nabla,i}(z)}-2h^\vee (A(z))^\nabla(z)-kh^\vee(A(z),A(z))
    \bigr),
    \end{split}
\]
where $(A(z))^\nabla(z)$ denotes $\sum_{k}A^\nabla_k(z)z^{-k-1}=\sum_{k,l}(A_k)^\nabla_lz^{-k-l-2}$. Now, let $x(z)\ceq x^\nabla(z)+(A(z),x)k$. Then, $x(z)$ satisfying the relation of $\wh{\frg}$ by \cref{prp:isom aff Lie alg}. Since in this notation, we write 
\[
    \begin{split}
    &\nop{x_i^\nabla(z)x^{\nabla,i}(z)}=\nop{x_i(z)x^i(z)}-2k(A(z))(z)+k^2(A(z),A(z)),\\
    &(A(z))^\nabla(z)=(A(z))(z)-(A(z),A(z))k
    \end{split}
\]
we have
\[
    L^\nabla_\mathrm{sug}(z)=\sum_i\frac{\nop{x_i(z)x^i(z)}}{2(k+h^\vee)}-(A(z))(z)+\frac{kh}{2}(A(z),A(z)).
\]

\begin{prp}\label{prp:conf str on tw aff}
    On $\nabla$-twisted $V^k(\frg)$-module $M$, the operator
    \[
        L^\nabla_\mathrm{sug}(z)=\sum_i\frac{\nop{x_i(z)x^i(z)}}{2(k+h^\vee)}-(A(z))(z)+\frac{kh}{2}(A(z),A(z)).
    \]
    gives a Virasoro action at central charge $c=k\dim(\frg)/(k+h^\vee)$. Here, $x(z)$ is an operator given by $x^\nabla(z)+(A(z),x)k$, satisfying the relation with $\wh{\frg}$. 
\end{prp}

\section*{Acknowledgments}
    The author would like to express his sincere gratitude to his advisor, S.\ Yanagida, for valuable discussions, helpful advice, and continued encouragement. The author is also grateful to A.\ Shimojima, N.\ Tajima, T.\ Yano, and S.\ Yoshida for their interest in the author's ideas and for valuable discussions.

    The author used ChatGPT (GPT-6, OpenAI) for language editing and its assistance in developing the proof presented in Appendix~A. Building on the author's original idea of applying Tannakian formalism to construct a $G$-valued formal solution, discussions with ChatGPT helped develop an argument using the fiber functor of horizontal sections and the triviality of torsors over an algebraically closed field. The author takes full responsibility for the correctness of the resulting proof.

    This work was supported by JSPS KAKENHI Grant Number JP25KJ1393.

\appendix

\section{Existence of $G$-valued solutions}\label{s: Existence of G-valued solutions}
In this appendix, we establish the existence of a $G$-valued fundamental solution of a regular singular $\frg$-connection, where $G$ is an affine algebraic group over $\bbC$ with Lie algebra $\frg$. Although the existence of fundamental solutions is well known for linear ordinary differential equations, it is not immediate that a fundamental solution can be chosen to take values in $G$. To address this issue, we use the special case of Tannakian formalism.

Let $\Rep(G)$ be the category of representations of $G$ and $\omega\colon \Rep(G)\to \cVect$ be the forgetful functor. If a $\bbC$-linear functor $\eta\colon\Rep(G)\to \cVect$ is an exact, faithful tensor functor, then it is called a \emph{fiber functor}. For example, $\omega$ is a standard example of a fiber functor. For two $\bbC$-linear tensor functors $\eta_1,\eta_2\colon\Rep(G)\to\cVect$, we denote by $\Isom^\otimes(\eta_1,\eta_2)$ the set of natural isomorphisms $\eta_1\to\eta_2$ compatible with tensors and $\Aut^\otimes(\eta_1)=\Isom^\otimes(\eta_1,\eta_1)$.  

For a $\frg$-connection $\nabla=\pd_z+A(z)$, we define a functor $\eta_\nabla\colon\Rep(G)\to\cVect$ by
\[
    \eta_\nabla(V)\ceq(V\otimes\wt{\shK}_z)^{\nabla_V},\quad \nabla_V=\pd_z+\rho(A(z)),
\]
for $V=(V,\rho)\in\Rep(G)$, where $\wt{\shK}_z$ is a commutative algebra defined by
\[
    \wt{\shK}_z=\bigoplus_{\alpha\in\bbC/\bbZ,m\geq 0}\bbC\dpr{z}z^\alpha(\log z)^m
\]
with derivation $\pd_z$ defined by 
\[
    \pd_zz^\alpha=\alpha z^{\alpha-1},\quad \pd_z\log z=\frac{1}{z}.
\]
\begin{lem}\label{lem:fiber functor}
    If $\nabla$ is a regular singular $\frg$-connection, then the functor $\eta_\nabla$ is a fiber functor. 
\end{lem}

\begin{proof}
    By the fundamental theorem for linear ordinary differential equations, we have
    \[
        \dim_\bbC\eta_\nabla(V)=\dim_\bbC(V)
    \]
    for any $V\in\Rep(G)$. By comparing ranks, the natural map $\eta_\nabla(V)\otimes\wt{\shK}_z\to \omega(V)\otimes\wt{\shK}_z$ is an isomorphism. Since the functor $V\mto \omega(V)\otimes\wt{\shK}_z$ is faithful, so is $\eta_\nabla$. 

    Next, we check the exactness of $\eta_\nabla$. It suffices to prove the right exactness. Let $f\colon V\to W$ be a surjective $G$-module homomorphism and $C$ be a cokernel of $\eta_\nabla(f)\colon \eta_\nabla(V)\to\eta_\nabla(W)$. By the following commutative diagram:
    \[
        \begin{tikzcd}
            \eta_\nabla(V)\otimes\wt{\shK}_z \ar[r,"\eta_\nabla(f)\otimes 1"]\ar[d,"\sim"' sloped] &\eta_\nabla(W)\otimes \wt{\shK}_z \ar[d,"\sim" sloped] \\
            V\otimes\wt{\shK}_z \ar[r,"f\otimes 1"'] &W\otimes\wt{\shK}_z, 
        \end{tikzcd}
    \]
    we have $C\otimes\wt{\shK}_z=0$, which implies $C=0$.

    Finally, we show that $\eta_\nabla$ is a tensor functor. Note that for $V,W\in\Rep(G)$, $\nabla_{V\otimes W}=\nabla_V\otimes\id_W+\id_V\otimes\nabla_W$. Thus, we obtain a natural injective map $\eta_\nabla(V)\otimes\eta_\nabla(W)\hookrightarrow\eta_\nabla(V\otimes W)$. By dimension comparison, this map is isomorphic.
\end{proof}

\begin{thm}\label{thm:existence of G-valued solutions}
    For any regular singular $\frg$-connection $\nabla$, there exists $F(z)\in G(\wt{\shK}_z)$ such that
    \[
        F(z)^{-1}\nabla F(z)=d.
    \]
\end{thm}

\begin{proof}
    By \cite[Theorem~1.2]{Wib25}, we have $\Isom^\otimes(\omega,\eta_\nabla)\neq\varnothing$. Since $\eta_\nabla(V)$ is canonically contained in $\omega(V)\otimes\wt{\shK}_z$ for any $V\in\Rep(G)$, we have
    \[
        \varnothing\neq\Isom^\otimes(\omega,\eta_\nabla)\subset\Isom^\otimes(\omega,\omega\otimes\wt{\shK}_z)\simeq \Aut^\otimes_{\wt{\shK}_z}(\omega\otimes\wt{\shK}_z).
    \]
    By \cite[Proposition~2.8]{DM82}, we have $G(\wt{\shK}_z)= \Aut^\otimes_{\wt{\shK}_z}(\omega_{\wt{\shK}_z})$. Thus, we obtain an element $F(z)\in G(\wt{\shK}_z)$ such that
    \[
        \rho(F(z))\colon V\lto (V\otimes\wt{\shK}_z)^\nabla  
    \]
    for any $V=(V,\rho)\in\Rep(G)$. This means that $\rho(F(z))$ is a fundamental solution. Taking a faithful representation of $G$, we conclude that $F(z)$ is the desired element.
\end{proof}


\bibliographystyle{mybstwithlabels}
\bibliography{references}
\end{document}

%% file: macro.tex
\usepackage{fullpage}
\usepackage[foot]{amsaddr}
\usepackage{mathtools,amssymb,mathabx,amsthm}
\usepackage{mathrsfs,euscript}
\usepackage{braket,slashed}
\usepackage[only llbracket,rrbracket,longmapsfrom]{stmaryrd}
\usepackage{tikz}
\usetikzlibrary{cd}
\usepackage{showkeys,comment}
\usepackage[l2tabu,orthodox]{nag}
\usepackage[all, warning]{onlyamsmath}

\usepackage[bookmarks=false,draft=false,breaklinks,colorlinks]{hyperref}

\usepackage{color}

\numberwithin{equation}{section}
\allowdisplaybreaks

{\par \vspace{\baselineskip}%
 \noindent \textbf{Acknowledgements.}}%
{\par \vspace{\baselineskip}}

\usepackage{enumitem}
\setenumerate{label=(\arabic*),nosep}
\setitemize{nosep}
\newlist{clist}{enumerate}{1}
\setlist*[clist]{label=(\roman*), nosep}

\usepackage{cleveref}

\crefname{thm}{Theorem}{Theorems}
\crefname{dfn}{Definition}{Definitions}
\crefname{prp}{Proposition}{Propositions}
\crefname{lem}{Lemma}{Lemmas}
\crefname{cor}{Corollary}{Corollaries}
\crefname{clm}{Claim}{Claims}
\crefname{fct}{Fact}{Facts}
\crefname{rmk}{Remark}{Remarks}
\crefname{eg}{Example}{Examples}
\crefname{figure}{Figure}{Figures}
\crefname{table}{Table}{Tables}
\crefname{section}{\S\!}{\S\S\!}
\crefname{subsection}{\S\!}{\S\S\!}
\crefname{subsubsection}{\S\!}{\S\S\!}
\crefname{appendix}{Appendix}{Appendices}
\crefname{equation}{}{}
\crefname{mainthm}{Theorem}{Theorems}

\theoremstyle{definition}
\newtheorem{thm}{Theorem}[section]
\newtheorem{dfn}[thm]{Definition}
\newtheorem{prp}[thm]{Proposition}
\newtheorem{lem}[thm]{Lemma}
\newtheorem{cor}[thm]{Corollary}

\newtheorem{rmk}[thm]{Remark}

\newtheorem*{rmk*}{Remark}
\newtheorem{mainthm}{Theorem}

\AddToHook{env/dfn/begin}{\crefalias{thm}{dfn}}
\AddToHook{env/prp/begin}{\crefalias{thm}{prp}}
\AddToHook{env/lem/begin}{\crefalias{thm}{lem}}
\AddToHook{env/cor/begin}{\crefalias{thm}{cor}}
\AddToHook{env/clm/begin}{\crefalias{thm}{clm}}
\AddToHook{env/fct/begin}{\crefalias{thm}{fct}}
\AddToHook{env/rmk/begin}{\crefalias{thm}{rmk}}
\AddToHook{env/eg/begin}{\crefalias{thm}{eg}}

\newcommand{\pd}{\partial}

\newcommand{\ul}{\underline}
\newcommand{\wh}{\widehat}
\newcommand{\wt}{\widetilde}
\newcommand{\ceq}{\coloneqq} 

\newcommand{\lto}{\longrightarrow}

\newcommand{\mto}{\mapsto}
\newcommand{\lmto}{\longmapsto}

\newcommand{\bbC}{\mathbb{C}}

\newcommand{\bbZ}{\mathbb{Z}}

\newcommand{\shK}{\EuScript{K}}

\newcommand{\shO}{\EuScript{O}}

\newcommand{\frg}{\mathfrak{g}}

\newcommand{\frsl}{\mathfrak{sl}}

\newcommand{\one}{\mathbf{1}}

\newcommand{\cMod}{\mathsf{Mod}}

\newcommand{\Rep}{\mathsf{Rep}}
\newcommand{\cVect}{\mathsf{Vect}}

\newcommand{\abs}[1]{\left| #1 \right|}
\newcommand{\dbr}[1]{\llbracket #1 \rrbracket} 
\newcommand{\dpr}[1]{(\!( #1 )\!)}

\newcommand{\nop}[1]{\,\substack{\scriptscriptstyle\circ \\ \scriptscriptstyle\circ} #1 \substack{\scriptscriptstyle\circ \\ \scriptscriptstyle\circ\,}}

\DeclareMathOperator{\id}{id}
\DeclareMathOperator{\GL}{GL}
\DeclareMathOperator{\SL}{SL}

\DeclareMathOperator{\Aut}{Aut}

\DeclareMathOperator{\Hom}{Hom}
\DeclareMathOperator{\Isom}{Isom}

\DeclareMathOperator{\LFie}{LFie}

\DeclareMathOperator{\shBer}{\EuScript{B}\kern-.1em\mathit{er}}
\DeclareMathOperator{\shDer}{\EuScript{D}\kern-.1em\mathit{er}}
\DeclareMathOperator{\shEnd}{\EuScript{E}\kern-.1em\mathit{nd}}
\DeclareMathOperator{\shHom}{\EuScript{H}\kern-.15em\mathit{om}}
\DeclareMathOperator{\shSym}{\EuScript{S}\kern-.1em\mathit{ym}}